\documentclass[12pt]{article}
\usepackage{amsmath,amsthm,amsfonts,amssymb,epic,eepic}
\makeatletter
\def\section{\@startsection {section}{1}{\z@}{-3.5ex plus -1ex minus
 -.2ex}{2.3ex plus .2ex}{\normalsize\bf}}
\def\@maketitle{\newpage
 \null
 \vskip 2em \begin{center}
 {\large \@title \par} \vskip 1.5em {\normalsize \lineskip .5em
\begin{tabular}[t]{c}\@author
 \end{tabular}\par}
 \vskip 1em {\normalsize \@date} \end{center}
 \par
 \vskip 1.5em}
\makeatother

\title{{\bf Generalized Longo-Rehren subfactors and $\alpha$-induction}}
\author{
{\sc Yasuyuki Kawahigashi}\\
Department of Mathematical Sciences\\
University of Tokyo, Komaba, Tokyo, 153-8914, JAPAN\\
e-mail: {\tt yasuyuki@ms.u-tokyo.ac.jp}}

\begin{document}
\maketitle

\def\isom{\cong}

\def\Ad{{\mathrm{Ad}}}
\def\Aut{{\mathrm{Aut}}}
\def\End{{\mathrm{End}}}
\def\Hom{{\mathrm{Hom}}}
\def\id{{\mathrm{id}}}
\def\Mor{{\hbox{Mor}}}
\def\NXN{{}_N{\mathcal{X}}_N}
\def\NXM{{}_N{\mathcal{X}}_M}
\def\MXN{{}_M{\mathcal{X}}_N}
\def\MXM{{}_M{\mathcal{X}}_M}
\def\o{{\mathrm{opp}}}
\def\op{{\mathrm{opp}}}
\def\opp{{\mathrm{opp}}}
\def\ti{\tilde}

\def\C{{\bf C}}
\def\N{{\bf N}}
\def\Q{{\bf Q}}
\def\R{{\bf R}}
\def\Z{{\bf Z}}

\def\a{\alpha}
\def\be{\beta}
\def\de{\delta}
\def\e{\varepsilon}
\def\Ep{\mathcal{E}}
\def\epsilon{\varepsilon}
\def\ga{\gamma}
\def\la{\lambda}
\def\si{\sigma}
\def\th{\theta}
\def\Ga{\Gamma}
\def\Th{\Theta}

\newtheorem{theorem}{Theorem}[section]
\newtheorem{lemma}[theorem]{Lemma}
\newtheorem{conjecture}[theorem]{Conjecture}
\newtheorem{corollary}[theorem]{Corollary}
\newtheorem{definition}[theorem]{Definition}
\newtheorem{assumption}[theorem]{Assumption}
\newtheorem{proposition}[theorem]{Proposition}
\newtheorem{remark}[theorem]{Remark}
\newtheorem{example}[theorem]{Example}

\def\emptyset{\varnothing}
\def\setminus{\smallsetminus}

\def\qed{{\unskip\nobreak\hfil\penalty50
\hskip2em\hbox{}\nobreak\hfil Q.E.D.
\parfillskip=0pt \finalhyphendemerits=0\par}\medskip}
\def\proof{\trivlist \item[\hskip \labelsep{\it Proof.\ }]}
\def\endproof{\null\hfill\qed\endtrivlist}

\begin{abstract}
We study the recent construction of subfactors by Rehren
which generalizes the Longo-Rehren subfactors.  We prove
that if we apply this construction to a non-degenerately
braided subfactor $N\subset M$ and $\a^\pm$-induction, then
the resulting subfactor is dual to the Longo-Rehren subfactor
$M\otimes M^\op \subset R$ arising from the entire system
of irreducible
endomorphisms of $M$ resulting from $\a^\pm$-induction.  As
a corollary, we solve a problem on existence of braiding
raised by Rehren negatively.
Furthermore, we generalize our previous study with Longo and M\"uger
on multi-interval subfactors arising from a completely rational
conformal net of factors on $S^1$ 
to a net of subfactors and show that the (generalized) Longo-Rehren
subfactors and $\a$-induction naturally appear in this context.
\end{abstract}

\section{Introduction}

In subfactor theory initiated by V. F. R. Jones \cite{J},
Ocneanu's construction of
asymptotic inclusions \cite{O1} have been studied by
several people as a subfactor analogue of the quantum double
construction.  (See \cite[Chapter 12]{EK} on general theory
of asymptotic inclusions.)  Popa's construction of symmetric
enveloping inclusions \cite{P} gives its generalizations and is
important in the analytic aspects of subfactor theory.  Longo and
Rehren gave another construction of subfactors in \cite{LR}
in the setting of sector theory \cite{L1,L2}
and Masuda \cite{M} has proved that the asymptotic inclusion
and the Longo-Rehren subfactor are essentially the same
constructions.  Izumi \cite{I1,I2} gave very detailed and 
interesting studies of the Longo-Rehren subfactors.
Recently, Rehren \cite{R2} gave a construction generalizing the
Longo-Rehren subfactor and we call the resulting subfactor a
{\sl generalized Longo-Rehren subfactor}.  This construction
uses certain extensions of systems of endomorphisms from
subfactors (of type III) to larger factors.  We will analyze
this construction in detail in this paper.
(This construction will be explained in more detail
in Section \ref{GLR} below.)

Longo and Rehren also defined such an extension of endomorphisms
for nets of subfactors in the same paper \cite[Proposition 3.9]{LR},
based on an old suggestion of Roberts
\cite{Ro}.  The essentially same construction of new endomorphisms
was also given in Xu \cite[page 372]{X1} and several very interesting
properties and examples were found
by him in \cite{X1,X2}.  We call this extension of endomorphisms
{\sl $\a$-induction}.   In this paper, we study the generalized
Longo-Rehren subfactors arising from $\a$-induction based on
the above works, 
B\"ockenhauer-Evans \cite{BE} and
our previous work \cite{BEK1,BEK2,BEK3}.
In the papers of Longo, Rehren, and Xu, they study nets of
subfactors and have a certain condition
arising from locality of the larger net, now called
{\sl chiral locality} as in \cite[Section 3.3]{BEK1},
but we do not assume this condition
in this paper.  We assume only a non-degenerate braiding 
in the sense of \cite{R1}.  (See \cite[Section 3.3]{BEK1}
for more on this matter.  We only need a braiding in order to define
$\a$-induction, but we also assume non-degeneracy in
this paper.  If we
start with a completely rational net on the circle in the
sense of \cite{KLM}, non-degeneracy of the braiding holds
automatically by \cite{KLM}.)  Izumi's work \cite{I1,I2}
on a half-braiding is closely related to theory of $\a$-induction
and a theory of induction for bimodules generalizing
these works has been recently given by Kawamuro \cite{KK}.

Results in \cite{BEK3} suggest that if we apply the
construction of the generalized Longo-Rehren subfactor
to $\a^\pm$-induction for $N\subset M$, then the resulting
subfactor $N\otimes N^\op\subset P$ would be dual to the
Longo-Rehren subfactor $M\otimes M^\op\subset R$ applied
to the system of endomorphisms of $M$ arising from
$\a^\pm$-induction.  In this paper we will prove that this
is indeed the case.  The proof involves several calculations
of certain intertwiners related to a half-braiding in the
sense of Izumi \cite{I1} arising from a relative braiding
in B\"ockenhauer-Evans \cite{BE}.  As
an application, we solve a problem on existence of braiding
raised by Rehren \cite{R2} negatively.

Furthermore, we generalize our previous study with Longo and M\"uger
\cite{KLM} on multi-interval subfactors arising from a completely rational
conformal net of factors on $S^1$ 
to a net of subfactors.  That is, we have studied ``multi-interval
subfactors'' arising from such a net
on $S^1$, whose definitions will be explained below,
and proved that the resulting subfactor is isomorphic to the Longo-Rehren
subfactor arising from all superselection sectors of the net
in \cite{KLM}.  We apply the construction of multi-interval
subfactors to conformal nets of {\sl subfactors} with finite index and
prove that the resulting subfactor is isomorphic to the Longo-Rehren
subfactor arising from the system of $\a$-induced endomorphisms.
We then also explain a relation of this result to the generalized
Longo-Rehren subfactors.

The results in Section \ref{GLR} were announced in \cite{K}.

\section{Generalized Longo-Rehren subfactors}
\label{GLR}

Let $N\subset M$ be a type III subfactor with finite index
and finite depth.  Let $\NXN,\NXM,\MXN,\MXM$ be finite systems
of irreducible morphisms of type $N$-$N$, $N$-$M$, $M$-$N$, $M$-$M$,
respectively and suppose that the four systems together make a
closed system under conjugations, compositions and irreducible
decompositions, and the inclusion map from $N$ into $M$
decomposes into irreducible $N$-$M$ morphisms within
$\NXM$,  as in \cite[Assumption 4.1]{BEK1}.
We assume that the system $\NXN$ is non-degenerately braided
as in \cite{R1},
\cite[Definition 2.3]{BEK1}.  Then we have positive and negative
$\a$-inductions, corresponding to positive and negative
braidings, and the system $\MXM$ is generated by
the both $\a$-inductions because of the
non-degeneracy as in \cite[Theorem 5.10]{BEK1}.  We do {\sl not}
assume the chiral locality condition, which arises from locality
of the larger net of factors, in this paper.  (See
\cite[Section 5]{BEK2} for more on the role of chiral locality.)

Now recall a new construction of subfactors due to Rehren \cite{R2}
arising from two systems of endomorphisms and two extensions
to the same factor as follows.

Let $\Delta$ be a system of endomorphisms of a type III factor
$N$ and consider a subfactor $N\subset M$ with finite index.
An {\sl extension} of $\Delta$ is a pair 
$(\iota,\a)$ where $\iota$ is the embedding map of $N$ into
$M$ and $\a$ is a map $\Delta\to\End(M)$, $\la\mapsto\a_\la$
satisfying the following properties.
\begin{enumerate}
\item Each $\a_\la$ has a finite dimension.
\item We have $\iota\la=\a_\la\iota$ for $\la \in\Delta$.
\item We have $\iota(\Hom(\la\mu,\nu))
\subset \Hom(\a_\la\a_\mu,\a_\nu)$ for $\la,\mu,\nu\in\Delta$.
\end{enumerate}

Next let $N_1$, $N_2$ be two subfactors  of a type III factor
$M$, $(\iota_1,\a^1)$ and $(\iota_2,\a^2)$ be two
extensions of finite systems $\Delta_1, \Delta_2$ of endomorphisms
of $N_1, N_2$ to $M$, respectively.  For $\la\in\Delta_1$ and
$\mu\in\Delta_2$, we set $Z_{\la,\mu}=\dim\Hom(\a^1_\la,\a^2_\mu)$.
Then Rehren proved in \cite{R2} that we have a subfactor 
$N_1\otimes N_2^\o\subset R$ such that the canonical endomorphism
restricted on $N_1\otimes N_2^\opp$ has a decomposition
$\bigoplus_{\la\in\Delta_1,\mu\in\Delta_2} Z_{\la,\mu}
\la\otimes\mu^\opp$ by constructing the corresponding $Q$-system
explicitly.  This is a generalization of the Longo-Rehren
construction \cite[Proposition 4.10]{LR}
in the sense that if $N_1=N_2=M$, Rehren's
$Q$-system coincides with the one given in \cite{LR}.
We call it a generalized Longo-Rehren subfactor.
The most natural example of such extensions seems to be the
$\a$-induction, and then
we can take $\Delta=\NXN$, $\a^1=\a^+$, $\a^2=\a^-$
for $\a$-induction from $N$ to $M$ based on a braiding
$\e^\pm$ on the system $\NXN$ and then
$Z_{\la,\mu}$ is the ``modular invariant'' matrix as
in \cite[Definition 5.5, Theorem 5.7]{BEK1}.

Our aim is to study the generalized Longo-Rehren
subfactor arising from $\NXN$ and
$\a^\pm$-induction in this way.
The result in \cite[Corollary 3.11]{BEK3} suggests that
this subfactor is dual to the Longo-Rehren subfactor arising
from $\MXM$, and we prove this is indeed the case.  For this
purpose, we study the Longo-Rehren subfactor arising from $\MXM$
first as follows.

Let $M\otimes M^\o\subset R$ be the Longo-Rehren subfactor
\cite[Proposition 4.10]{LR}
arising from the system $\MXM$ on $M$ and
$(\Gamma,V,W)$ be the corresponding $Q$-system \cite{L3}.
(Actually, the subfactor we deal with here is the dual to
the original one constructed in \cite[Proposition 4.10]{LR}.  This dual
version is called the Longo-Rehren subfactor in
\cite{BEK3}, \cite{KLM}.)
That is, we have that $\Gamma\in \End(R)$ is the canonical endomorphism
of the subfactor, $V\in\Hom(\id,\Gamma)\subset R$, and
$W\in\Hom(\Gamma,\Gamma^2)$.  We also have
\begin{eqnarray*}
W&\in& M\otimes M^\o,\\
R&=&(M\otimes M^\o)V,\\
W^*V&=&\Gamma(V^*)W=w^{-1/2},\\
W^*\Gamma(W)&=&WW^*,\\
\Gamma(W)W&=&W^2,
\end{eqnarray*}
where $w=\sum_{\be\in\MXM} d^2_\be$ is
the global index of the system $\MXM$ and equal to the
index $[R:M\otimes M^\o]$.
By the definition of the original Longo-Rehren subfactor
in \cite{LR}, the $Q$-system $(\Theta,W,\Ga(V))$ is given
as follows.  We have
$$\Th(x)=\sum_{\be\in\MXM} W_\be (\be\otimes\be^\o(x)) W_\be^*,
\quad\mathrm{for \ }x \in M\otimes M^\o,$$
where $\Th$ is the dual canonical endomorphisms, the restriction
of $\Ga$ to $M\otimes M$,
the family $\{W_\be\}$ is that of isometries with mutually
orthogonal ranges satisfying $\sum_{\be\in\MXM} W_\be W_\be^*=1$,
and also have 
\begin{eqnarray}
\label{LRQ}
\Ga(V)&=&\sum_{\be_1,\be_2,\be_3\in\MXM}\sqrt{\frac{d_1 d_2}{w d_3}}
\Ga(W_{\be_2})W_{\be_1} T_{\be_1\be_2}^{\be_3} W^*_{\be_3},\\
T_{\be_1\be_2}^{\be_3} &=& \sum_{l=1}^{N_{12}^3} 
T_{\be_1\be_2, l}^{\be_3} \otimes j(T_{\be_1\be_2, k}^{\be_3})
\in M\otimes M^\o,
\end{eqnarray}
by definition of the Longo-Rehren subfactor \cite{LR}, where
$\{T_{\be_1\be_2, l}^{\be_3}\}_l$ is an orthogonal basis in
$\Hom(\be_3,\be_1\be_2)\subset M$, $N_{ij}^k$ is the
structure constant $\dim\Hom(\be_k,\be_i\be_j)$, 
$d_j=d_{\be_j}$ is the statistical dimension of $\be_j$, and
$j$ is the anti-isomorphism $x\in M\mapsto x^*\in M^\o$.
Starting from this explicit expression of the $Q$-system
$(\Theta,W,\Ga(V))$, we would like to write down the
$Q$-system $(\Ga,V, W)$ explicitly and identify it with the
$Q$-system given by the construction of Rehren \cite{R2}.

First, by \cite[Theorem 3.9]{BEK3}, we know that
$$[\Ga]=\bigoplus_{\la_1,\la_2\in\NXN}Z_{\la_1,\la_2}
[\eta(\a_{\la_1}^+,+) \eta^\o(\a_{\la_2}^-,-)],$$
where $[\hphantom{X}]$ represents the sector class of an endomorphism,
$\eta(\hphantom{X},\hphantom{X})$ is the extension of an endomorphism of
$M\otimes M^\o$ to $R$ with a half-braiding by Izumi \cite{I1},
$\a^\pm$ is the $\a$-induction,
the notations here follow those of \cite{BEK3}, and 
$Z_{\la_1\la_2}=\dim\Hom(\a_{\la_1}^+,\a_{\la_2}^-)$
is the ``modular invariant''
as in \cite[Definition 5.5]{BEK1}.  (Recall that we now assume
non-degeneracy of the braiding on $\NXN$.)
Furthermore, by \cite[Corollary 3.10]{BEK3}, we have equivalence
of two $C^*$-tensor categories of
$\{\eta(\a_{\la}^+,+) \eta^\o(\a_{\mu}^-,-)\}$ on $R$
and $\{\la\otimes\mu^\o\}$ on $N\otimes N^\o$, thus the
canonical endomorphisms of the two $Q$-systems are naturally
identified.  So we will next compute $V, W$ explicitly and identify
them with the intertwiners in Rehren's $Q$-system.
(Note that it does not matter that two von Neumann algebras $R$
and $N\otimes N^\o$ are different, since only the equivalence class
of $C^*$-tensor categories matters in the construction of the
(generalized) Longo-Rehren subfactors.)

We next closely follow Izumi's arguments in \cite[Section 7]{I1}.
First we have the following lemma.

\begin{lemma}\label{Le1}
For an operator $X\in M\otimes M^\o$, $XV\in R$ is in
$$\Hom(\eta(\a_{\la_1}^+,+)\eta^\o(\a_{\la_2}^-,-),\Ga)$$
if and only if we have the following two conditions.
\begin{enumerate}
\item $X\in\Hom(\Th(\a_{\la_1}^+\otimes\a_{\la_2}^{-,\o}),\Th)$.
\item $X\Ga(U^*)\Ga(V)=\Ga(V)X$, where
$$U=\sum_{\be\in\MXM} W_\be(\Ep^+_{\la_1}(\be)
\otimes j(\Ep_{\la_2}^-(\be)))
(\a_{\la_1}^+\otimes\a_{\la_2}^{-,\o})(W_\be^*),$$
and $\Ep^\pm$ is the half-braiding defined in
\cite[Section 3]{BEK3}.
\end{enumerate}
\end{lemma}

\begin{proof}
By a standard argument similar to the one in the proof of
\cite[Proposition 7.3]{I1}, we easily get the conclusion.
\end{proof}

Next, we rewrite the second condition in the above Lemma as
follows.  Using the definition of $\Ga(V)$ as in (\ref{LRQ}),
we have
\begin{eqnarray*}
&&\sum_{\be_1,\be_2,\be_3\in\MXM}X\Ga(U^*)\sqrt{\frac{d_1 d_2}{w d_3}}
\Ga(W_{\be_2})W_{\be_1} T_{\be_1\be_2}^{\be_3} W^*_{\be_3}\\
&=& \sum_{\be_4,\be_5,\be_6\in\MXM}\sqrt{\frac{d_4 d_5}{w d_6}}
\Ga(W_{\be_5})W_{\be_4} T_{\be_4\be_5}^{\be_6} W^*_{\be_6}X,
\end{eqnarray*}
which is equivalent to the following equations for all $\be_3,
\be_4,\be_5\in\MXM$.
\begin{eqnarray*}
&&\sum_{\be_1,\be_2\in\MXM}W_{\be_4}^*\Th(W_{\be_5})X\Th(U^*)
\sqrt{\frac{d_1 d_2}{d_3}}
\Th(W_{\be_2})W_{\be_1} T_{\be_1\be_2}^{\be_3}\\
&=& \sum_{\be_6\in\MXM}\sqrt{\frac{d_4 d_5}{d_6}}
T_{\be_4\be_5}^{\be_6} W^*_{\be_6}X W_{\be_3}.
\end{eqnarray*}
Assuming the first condition in Lemma \ref{Le1},
we compute the left hand side of this equation as follows.
\begin{eqnarray*}
&&\sum_{\be_1,\be_2\in\MXM}
\sqrt{\frac{d_1 d_2}{d_3}}
W_{\be_4}^* X \Th((\a_{\la_1}^+\otimes\a_{\la_2}^{-,\o})
(W^*_{\be_5}) U^*
W_{\be_2})W_{\be_1} T_{\be_1\be_2}^{\be_3}\\
&=&\sum_{\be_1,\be_2\in\MXM}
\sqrt{\frac{d_1 d_2}{d_3}}
W_{\be_4}^* X W_{\be_1}(\be_1\otimes\be_1^\o)
((\a_{\la_1}^+\otimes\a_{\la_2}^{-,\o})(W^*_{\be_5}) U^*
W_{\be_2}) T_{\be_1\be_2}^{\be_3}\\
&=&\sum_{\be_1,\be_2\in\MXM}
\sqrt{\frac{d_1 d_2}{d_3}}
W_{\be_4}^* X W_{\be_1}(\be_1\otimes\be_1^\o)
(\sum_{\be\in\MXM}(\a_{\la_1}^+\otimes\a_{\la_2}^{-,\o})(W^*_{\be_5})
(\a_{\la_1}^+\otimes\a_{\la_2}^{-,\o})\\
&&\qquad\qquad\qquad\times
(W_{\be})(\Ep^+_{\la_1}(\be)^*\otimes 
j(\Ep_{\la_2}^-(\be))^*)W_\be^* W_{\be_2})
T_{\be_1\be_2}^{\be_3}\\
&=&\sum_{\be_1\in\MXM}
\sqrt{\frac{d_1 d_5}{d_3}}
W_{\be_4}^* X W_{\be_1}(\be_1\otimes\be_1^\o)
(\Ep^+_{\la_1}(\be_5)^*\otimes j(\Ep_{\la_2}^-(\be_5))^*)
T_{\be_1\be_5}^{\be_3}.
\end{eqnarray*}
That is, our equation is now
\begin{eqnarray}
\label{Eq1}
&&\sum_{\be_1\in\MXM}
\sqrt{\frac{d_1}{d_3}}
W_{\be_4}^* X W_{\be_1}(\be_1\otimes\be_1^\o)
(\Ep^+_{\la_1}(\be_5)^*\otimes j(\Ep_{\la_2}^-(\be_5))^*)
T_{\be_1\be_5}^{\be_3}\nonumber\\
&=& \sum_{\be_6\in\MXM}\sqrt{\frac{d_4}{d_6}}
T_{\be_4\be_5}^{\be_6} W^*_{\be_6}X W_{\be_3}
\end{eqnarray}
for all $\be_3, \be_4,\be_5\in\MXM$.
Now set $\be_4=\id$ in this equation.  Then on the left
hand side, we have a term $W_0^* X W_{\be_1}$, which is in
$\Hom((\be_1\otimes \be_1^\o)(\a_{\la_1}^+
\otimes\a_{\la_2}^{-,\o}),
\id_{M\otimes M^\o})$.  

Now setting $X_{\be_1}=W_0^* X W_{\be_1}$, we get
$$\sum_{\be_1\in\MXM}
\sqrt{\frac{d_1}{d_3}}
X_{\be_1}(\be_1\otimes\be_1^\o)
(\Ep^+_{\la_1}(\be_5)^*\otimes j(\Ep_{\la_2}^-(\be_5))^*)
T_{\be_1\be_5}^{\be_3}
=\sqrt{\frac{1}{d_5}} W^*_{\be_5}X W_{\be_3}$$
for any $\be_3,\be_5\in\MXM$
from the equation (\ref{Eq1}), and this implies
\begin{equation}\label{Eq2}
X=\sum_{\be_1,\be_3,\be_5\in\MXM}
\sqrt{\frac{d_1 d_5}{d_3}} 
W_{\be_5}X_{\be_1}(\be_1\otimes\be_1^\o)
(\Ep^+_{\la_1}(\be_5)^*\otimes j(\Ep_{\la_2}^-(\be_5))^*)
T_{\be_1\be_5}^{\be_3} W^*_{\be_3}.
\end{equation}

Consider the linear map sending $X\in M\otimes M^\o$ with 
$$XV\in\Hom(\eta(\a_{\la_1}^+,+)
\eta^\o(\a_{\la_2}^-,-),\Ga)$$
to 
$$(W_0^* X W_{\be})_{\be}\in\bigoplus_{\be\in\MXM}
\Hom(\be\a_{\la_1}^+,\id)\otimes
\Hom(\be^\o \a_{\la_2}^{-,\o},\id).$$
The dimensions of the space of such $X$ and the space
$$\bigoplus_{\be\in\MXM}
\Hom(\be\a_{\la_1}^+,\id)\otimes
\Hom(\be\o \a_{\la_2}^{-,\o},\id)$$
are both equal to $Z_{\la_1\la_2}$, 
and this map is injective by the equation (\ref{Eq2}),
so this map is also surjective.  That is, a general form of
such an $X$ is determined now by the equation (\ref{Eq2}), where
$X_{\be}$'s are now arbitrary intertwiners in
$\Hom(\be\a_{\la_1}^+,\id)\otimes
\Hom(\be^\o \a_{\la_2}^{-,\o},\id)$.
Fix $\la_1,\la_2\in\NXN,\be\in\MXM$ and
$l_1$ and $l_2$ be indices in the set
$\{1,2,\dots,\dim\Hom(\be\a^+_{\la_1},\id)\}$,
$\{1,2,\dots,\dim\Hom(\be\a^-_{\la_2},\id)\}$ respectively.
Following \cite{R2}, we use the letter $l$ for the
multi-index $(\la_1,\la_2,\be,l_1,l_2)$.
Note that in order for us to get a non-trivial index, that is,
$l_1>0, l_2>0$, the endomorphism $\be$ must be {\sl ambichiral}
in the sense that it appears in irreducible decompositions of
both $\a^+$-induction and $\a^-$-induction as in \cite{BEK1}.
Let $\{T^+_{l_1}\}_{l_1}$ and
$\{T^-_{l_2}\}_{l_2}$ be orthonormal bases of
$\Hom(\a^+_{\la_1},\bar\be)$ 
and $\Hom(\a^-_{\la_2},\bar\be_2)$, respectively.

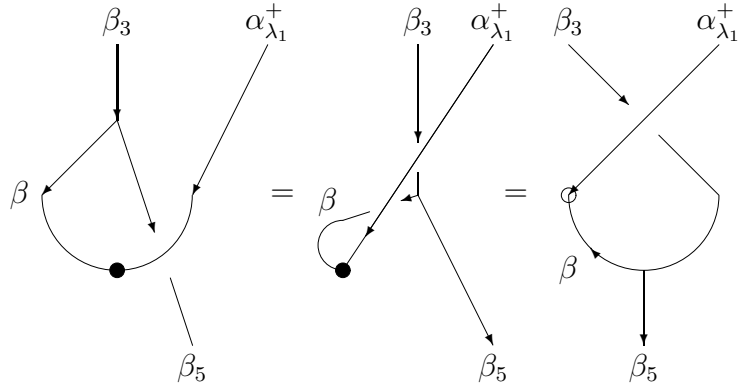
\begin{figure}[htb]
\begin{center}
\unitlength 1mm
\begin{picture}(110,60)
\thinlines
\put(20,50){\vector(0,-1){10}}
\put(20,40){\vector(-1,-1){10}}
\put(40,50){\vector(-1,-2){10}}
\put(20,40){\vector(1,-3){5}}
\put(30,10){\line(-1,3){3}}
\put(20,30){\arc{20}{0}{3.141593}}
\put(20,20){\circle*{2}}
\put(60,30){\vector(1,-2){10}}
\put(70,50){\vector(-2,-3){17}}
\put(70,50){\line(-2,-3){20}}
\put(60,50){\vector(0,-1){13}}
\put(60,30){\line(0,1){3}}
\put(60,30){\vector(-3,-1){2.2857}}
\put(50,26.667){\line(3,1){3.5714}}
\put(50,20){\circle*{2}}
\put(50,23.3333){\arc{6.6666}{1.570796}{4.712389}}
\put(90,20){\vector(0,-1){10}}
\put(82.928933,22.928933){\vector(-1,1){0}}
\put(90,30){\arc{20}{0}{3.141593}}
\put(80,30){\circle{2}}
\put(100,50){\vector(-1,-1){20}}
\put(80,50){\vector(1,-1){8}}
\put(92,38){\line(1,-1){8}}
\put(42,30){\makebox(0,0){$=$}}
\put(73,30){\makebox(0,0){$=$}}
\put(20,53){\makebox(0,0){$\be_3$}}
\put(40,53){\makebox(0,0){$\a^+_{\la_1}$}}
\put(60,53){\makebox(0,0){$\be_3$}}
\put(70,53){\makebox(0,0){$\a^+_{\la_1}$}}
\put(80,53){\makebox(0,0){$\be_3$}}
\put(100,53){\makebox(0,0){$\a^+_{\la_1}$}}
\put(30,7){\makebox(0,0){$\be_5$}}
\put(70,7){\makebox(0,0){$\be_5$}}
\put(90,7){\makebox(0,0){$\be_5$}}
\put(7,30){\makebox(0,0){$\be$}}
\put(48,29){\makebox(0,0){$\be$}}
\put(80,20){\makebox(0,0){$\be$}}
\end{picture}
\end{center}
\caption{An application of the braiding-fusion equation}
\label{Fg1}
\end{figure}

We now study some intertwiners using  a graphical calculus
in \cite[Section 3]{BEK1}.
First note that we have identities as in Fig.~\ref{Fg1}
by the braiding-fusion equation \cite[Definition 4.2]{I1},
\cite[Definition 2.2 2]{BEK3}
for a half-braiding, where crossings in the picture
represent the half-braidings and the black and white
small circles represent intertwiners in
$\Hom(\be\a^+_{\la_1},\id)$ and
$\Hom(\a^+_{\la_1},\bar\be)$ respectively.
(See \cite[Section 3]{BEK1} for interpretations of the graphical
calculus.  Here and below, a triple point, a black or while
small circle always represents an isometry or a co-isometry.
One has to be careful that we have a normalizing constant involving
the fourth roots of statistical dimensions as in
\cite[Figures 7,9]{BEK1}.  From now on, we drop orientations
of wires, which should causes no confusions.)
We also have the following lemma to relate
these two intertwiners.

\begin{figure}[htb]
\begin{center}
\unitlength 0.7mm
\begin{picture}(40,50)
\thinlines
\put(15,30){\arc{10}{3.141593}{6.283186}}
\put(15,20){\arc{10}{0}{3.141593}}
\put(30,10){\line(0,1){10}}
\put(30,30){\line(0,1){10}}
\put(10,20){\line(0,1){10}}
\put(20,20){\line(1,1){10}}
\put(20,30){\line(1,-1){4}}
\put(30,20){\line(-1,1){4}}
\put(10,25){\circle*{2}}
\put(30,44){\makebox(0,0){$\be$}}
\put(30,6){\makebox(0,0){$\a^+_\la$}}
\put(6,25){\makebox(0,0){$T_j^*$}}
\end{picture}
\end{center}
\caption{The intertwiner $\hat T_j$}
\label{Fg9}
\end{figure}
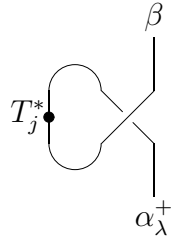

\begin{lemma}
\label{orth}
Let $T_j\in \Hom(\be,\a^+_\la)$ and define
$\hat T_j\in \Hom(\be,\a^+_\la)$
by the graphical expression in Fig.~\ref{Fg9}.
Then we have $T_k^* T_j=\hat T_k^* \hat T_j$.
\end{lemma}

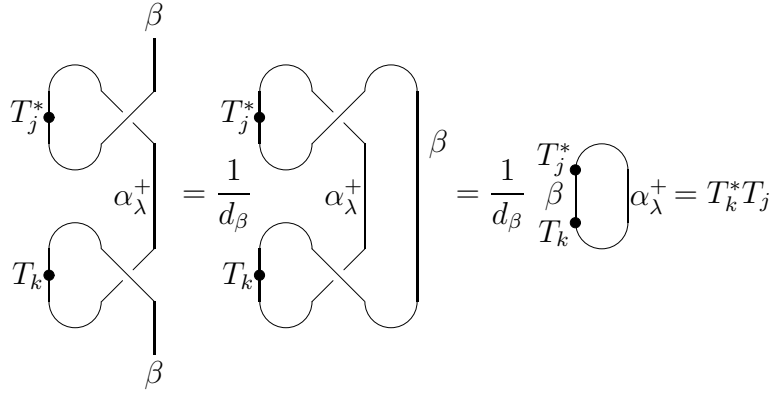
\begin{figure}[htb]
\begin{center}
\unitlength 0.7mm
\begin{picture}(130,80)
\thinlines
\put(15,20){\arc{10}{0}{3.141593}}
\put(15,50){\arc{10}{0}{3.141593}}
\put(75,20){\arc{10}{0}{3.141593}}
\put(55,20){\arc{10}{0}{3.141593}}
\put(55,50){\arc{10}{0}{3.141593}}
\put(115,35){\arc{10}{0}{3.141593}}
\put(15,30){\arc{10}{3.141593}{6.283186}}
\put(15,60){\arc{10}{3.141593}{6.283186}}
\put(55,30){\arc{10}{3.141593}{6.283186}}
\put(55,60){\arc{10}{3.141593}{6.283186}}
\put(75,60){\arc{10}{3.141593}{6.283186}}
\put(115,45){\arc{10}{3.141593}{6.283186}}
\put(30,10){\line(0,1){10}}
\put(30,60){\line(0,1){10}}
\put(30,30){\line(0,1){20}}
\put(10,20){\line(0,1){10}}
\put(10,50){\line(0,1){10}}
\put(70,30){\line(0,1){20}}
\put(50,20){\line(0,1){10}}
\put(50,50){\line(0,1){10}}
\put(80,20){\line(0,1){40}}
\put(110,35){\line(0,1){10}}
\put(120,35){\line(0,1){10}}
\put(20,50){\line(1,1){10}}
\put(20,60){\line(1,-1){4}}
\put(30,50){\line(-1,1){4}}
\put(60,50){\line(1,1){10}}
\put(60,60){\line(1,-1){4}}
\put(70,50){\line(-1,1){4}}
\put(20,30){\line(1,-1){10}}
\put(20,20){\line(1,1){4}}
\put(30,30){\line(-1,-1){4}}
\put(60,30){\line(1,-1){10}}
\put(60,20){\line(1,1){4}}
\put(70,30){\line(-1,-1){4}}
\put(10,25){\circle*{2}}
\put(10,55){\circle*{2}}
\put(50,25){\circle*{2}}
\put(50,55){\circle*{2}}
\put(110,35){\circle*{2}}
\put(110,45){\circle*{2}}
\put(30,74){\makebox(0,0){$\be$}}
\put(30,6){\makebox(0,0){$\be$}}
\put(84,50){\makebox(0,0){$\be$}}
\put(106,40){\makebox(0,0){$\be$}}
\put(26,40){\makebox(0,0){$\a^+_\la$}}
\put(66,40){\makebox(0,0){$\a^+_\la$}}
\put(124,40){\makebox(0,0){$\a^+_\la$}}
\put(6,55){\makebox(0,0){$T_j^*$}}
\put(6,25){\makebox(0,0){$T_k$}}
\put(46,55){\makebox(0,0){$T_j^*$}}
\put(46,25){\makebox(0,0){$T_k$}}
\put(106,47){\makebox(0,0){$T_j^*$}}
\put(106,33){\makebox(0,0){$T_k$}}
\put(42,40){\makebox(0,0)
{$=\displaystyle\frac{1}{d_\be}$}}
\put(94,40){\makebox(0,0)
{$=\displaystyle\frac{1}{d_\be}$}}
\put(138,40){\makebox(0,0){$=T_k^* T_j$}}
\end{picture}
\end{center}
\caption{The inner product $\hat T_k^* \hat T_j$}
\label{Fg10}
\end{figure}

\begin{proof}
We compute as in Fig.~\ref{Fg10}.
\end{proof}

Based on this, we set
$$S_{\be_1}^{\be_2\be_3} = \sum_{k=1}^{N_{23}^1} 
(T_{\be_2\be_3, k}^{\be_1})^* \otimes j(T_{\be_2\be_3, k}^{\be_1})^*
\in M\otimes M^\o$$ 
and we now define $X_l\in M\otimes M^\o$ as follows.
\begin{equation}
\label{Eq3}
X_l=\sqrt{d_{\la_1}d_{\la_2}}
\sum_{\be_3,\be_5\in\MXM}
\sqrt{\frac{d_3}{d_1 d_5}} W_{\be_5}
S^{\bar\be_1\be_3}_{\be_5}(T^+_{l_1}\otimes j(T^-_{l_2}))
(\Ep^+_{\la_1}(\be_3)^*\otimes j(\Ep_{\la_2}^-(\be_3))^*) W^*_{\be_3}.
\end{equation}
Then by the equation (\ref{Eq2}),
the operator $U_l\in R$ defined by $U_l=X_l V$ is in 
$$\Hom(\eta(\a_{\la_1}^+,+) \eta^\o(\a_{\la_2}^-,-),\Ga)$$
and $\{U_l\}_{\be,l_1,l_2}$ is a linear basis of
this intertwiner space.
We next prove that $\{U_l\}_{\be,l_1,l_2}$ is actually an
orthonormal basis with respect to the usual inner product.
Recall that for
$$s,t\in \Hom(\eta(\a_{\la_1}^+,+) \eta^\o(\a_{\la_2}^-,-),\Ga),$$
we have
$$E_{M\otimes M^\o}(st^*)=\frac{d_{\la_1}d_{\la_2}}{w}
t^*s\in\C$$
because $d_{\eta(\a_{\la_1}^+,+) \eta^\o(\a_{\la_2}^-,-)}=
d_{\la_1}d_{\la_2}$.  (See \cite[Lemma 3.1 (i)]{I1}.)
We then have
\begin{eqnarray*}
E_{M\otimes M^\o}(U_l U_{l'}^*)&=&
\frac{1}{w}X_l X^*_{l'}\\
&=&\de_{ll'}\frac{d_{\la_1}d_{\la_2}}{w}
\sum_{\be_3,\be_5\in\MXM}\frac{d_3}{d_1d_5} N_{13}^5
W_{\be_5}W^*_{\be_5}\\
&=& \de_{ll'}\frac{d_{\la_1}d_{\la_2}}{w},
\end{eqnarray*}
and this proves that $\{U_l\}_{\be,l_1,l_2}$ is indeed an
orthonormal basis. 
This also shows that we have
$$\phi_\Th(X_m^* X_l)=W^* E_{\Ga(R)}(X^*_m X_l)W=
W^*\Ga(U^*_m U_l)W=\de_{lm},$$
where $\phi_\Th$ is the standard left inverse of $\Th$.
(See \cite{LRo} for a general theory of left inverses.)

Let 
$l=(\la_1,\la_2,\be'_1,m_1,m_2)$,
$m=(\mu_1,\mu_2,\be''_1,m_1,m_2)$,
$n=(\nu_1,\nu_2,\be_1,n_1,n_2)$ be
multi-indices as above.  We compute
$E_{\Ga(R)}(X_m^* X_l^* X_n)$ as follows.
\begin{eqnarray*}
E_{\Ga(R)}(X_m^* X_l^* X_n)
&=& \Ga(V^* X_m^* X_l^* X_n V)\\
&=& \Ga(w^{1/2} V^* X_m^* \Ga(V^*)W X_l^* X_n V)\\
&=& \Ga(w^{1/2} V^* X_m^* \Ga(V^* X_l^*)W X_n V)\\
&=& \Ga(w^{1/2} U_m^* \Ga(U_l^*)W U_n).
\end{eqnarray*}
Based on this, we set
$$Y_{lm}^n=w^{-1/2}V^* X_m^* X^*_l X_n V=U_m^* \Ga(U_l^*)W U_n\in R$$
and then this is an element in
$$\Hom(\eta(\a_{\nu_1}^+,+) \eta^\o(\a_{\nu_2}^-,-),
\eta(\a_{\mu_1}^+,+) \eta^\o(\a_{\mu_2}^-,-)
\eta(\a_{\la_1}^+,+) \eta^\o(\a_{\la_2}^-,-)),$$
which is then contained in 
$$\Hom(\nu_1,\mu_1\la_1)\otimes\Hom(\nu_2,\mu_2\la_2)^\o
\subset N\otimes N^\o\subset M\otimes M^\o$$
by \cite[Theorem 3.9]{BEK3}.  That is, we now have
$$E_{\Ga(R)}(X_m^* X_l^* X_n)=
w^{1/2}\Th(Y_{lm}^n)\in\Th(M\otimes M^\o)$$ and
\begin{equation}
\label{Eq3+}
\phi_\Th(X_m^* X_l^* X_n)=V^* X_m^* X^*_l X_n V.
\end{equation}

\begin{proposition}
In the above setting, the $Q$-system $(\Gamma, V, W)$ is given as
follows.
\begin{eqnarray}
\label{Eq4}
\Ga(x)&=&\sum_l U_l 
(\eta(\a_{\la_1}^+,+) \eta^\o(\a_{\la_2}^-,-))(x) U_l^*,
\quad\mathrm{for\ }x\in R,\\
\label{Eq5}
V&=&U_{(0,0,0,1,1)},\\
\label{Eq6}
W&=&\sum_{l,m,n}\Ga(U_l) U_m Y^n_{lm} U_n^*.
\end{eqnarray}
\end{proposition}

\begin{proof}
Since $\{U_l\}_{\be_1,l_1,l_2}$ is an orthonormal basis of
$$\Hom(\eta(\a_{\la_1}^+,+) \eta^\o(\a_{\la_2}^-,-),\Ga),$$
we get the fist identity (\ref{Eq4}).
By the definition (\ref{Eq3}) of $X_l$, we have
$X_{(0,0,0,1,1)}=1$, hence $U_{(0,0,0,1,1)}=V$, which is (\ref{Eq5}).
Since $Y_{lm}^n=U_m^* \Ga(U_l^*)W U_n$, we get
(\ref{Eq6}).
\end{proof}

Next we further compute $Y_{lm}^n$.  We first have
\begin{eqnarray*}
Y_{lm}^n &=& W^* \Th(Y_{lm}^n) W\\
&=& w^{-1/2} W^* E_{\Ga(R)}(X_m^* X_l^* X_n) W\\
&=& w^{-1/2} \phi_\Th(X_m^* X_l^* X_n)\\
&=& \sum_{\be\in \MXM}\frac{d^2_\be}{w^{3/2}}
(\phi_\be\otimes\phi^\o_\be)
(W_\be^* X_m^* X_l^* X_n W_\be),
\end{eqnarray*}
where $\phi_\be$ is the standard left inverse of $\be$.
In this expression, we compute the term
$W_\be^* X_m^* X_l^* X_n W_\be$ as follows.
\begin{eqnarray*}
&& W_\be^* X_m^* X_l^* X_n W_\be\\
&=& \sqrt{d_{\la_1}d_{\la_2}d_{\mu_1}d_{\mu_2}d_{\nu_1}d_{\nu_2}}
\sum_{\be'_3,\be_5\in\MXM}
\frac{d_\be}{d_{\be_5}\sqrt{d_{\be''_1}d_{\be'_1}d_{\be_1}}}\\
&&\qquad\qquad\times
(\Ep^+_{\mu_1}(\be)\otimes j(\Ep_{\mu_2}^-(\be)))
((T^+_{m_1})^*\otimes j(T^-_{m_2})^*)
(S^{\bar\be''_1\be}_{\be'_3})^*\\
&&\qquad\qquad\times
(\Ep^+_{\la_1}(\be'_3)\otimes j(\Ep_{\la_2}^-(\be'_3)))
((T^+_{l_1})^*\otimes j(T^-_{l_2})^*)
(S^{\bar\be'_1\be'_3}_{\be_5})^*\\
&&\qquad\qquad\times
S^{\bar\be_1\be}_{\be_5}(T^+_{n_1}\otimes j(T^-_{n_2}))
(\Ep^+_{\nu_1}(\be)^*\otimes j(\Ep_{\nu_2}^-(\be))^*
\end{eqnarray*}

Our aim is to show that our $Y_{lm}^n$ coincides
with Rehren's ${\mathcal T}_{lm}^n$
in \cite[page 400]{R2}.  Our $Y_{lm}^n$ is already in
$\Hom(\nu_1,\mu_1\la_1)\otimes\Hom(\nu_2,\mu_2\la_2)^\o$ as
in Rehren's ${\mathcal T}_{lm}^n$.  So we expand
our $Y_{lm}^n$ with respect to the basis
$\{\tilde T_e = T^1_{e_1}\otimes j(T^2_{e_2})\}_{e=(e_1,e_2)}$,
where $\{T^1_{e_1}\}_{e_1}$, $\{T^2_{e_2}\}_{e_2}$ are bases for
$\Hom(\nu_1,\mu_1\la_1)$,
$\Hom(\nu_2,\mu_2\la_2)$, respectively.  We will prove that
the coefficients of $Y_{lm}^n$ for such an expansion
coincide with Rehren's coefficients $\zeta^n_{lm,e_1,e_2}$
in \cite[page 400]{R2}.

Let $S^+_l=S^+_{\be_1,\la_1,l_1}\in\Hom(\be_1,\a^+_{\la_1})$
be isometries so that $\{S^+_{\be_1,\la_1,l_1}\}_{l_1}$ gives
an orthonormal basis in $\Hom(\be_1,\a^+_{\la_1})$.  Similarly
we choose $S^-_l=S^-_{\be_1,\la_2,l_2}\in\Hom(\be_1,\a^-_{\la_2})$.
Rehren puts an inner product in $\Hom(\a^+_{\la_1},\a^-_{\la_2})$
in \cite[page 400]{R2}.   When we decompose this space as
$\bigoplus_{\beta\in\MXM^{0}} \Hom(\a^+_{\la_1},\beta)\otimes
\Hom(\beta,\a^-_{\la_2})$, Rehren's normalization implies that
his orthonormal basis consists
of intertwiners of the form 
$\sqrt{d_{\la_1}/d_{\beta}}S^-_l{S^+_l}^*$, where $S^\pm_l$
are isometries as above.
This implies that Rehren's $\zeta^n_{lm,e_1,e_2}$ is given as follows.
\begin{equation}
\label{Rc}
\sqrt{\frac{d_{\la_1}d_{\la_2}d_{\mu_1}d_{\mu_2}d_{\be''_1}}
{wd_{\nu_1}d_{\nu_2}d_{\be_1}d_{\be'_1}}}
S^{+*}_n (T^2_{e_1})^*((S_l^+ S_l^{-*})\times(S_m^+ S_m^{-*}))
T^1_{e_2} S_n^-.
\end{equation}

Note that we have
\begin{equation}
\label{Eq7}
E_{M\otimes M^\o}(X_n V \ti T_e V^* X_m^* X_l^*)=
\frac{d_{\nu_1} d_{\nu_2}}{w}\ti T_e V^* X_m^* X_l^* X_n V,
\end{equation}
where we have
$\tilde T_e = T^1_{e_1}\otimes j(T^2_{e_2})$ as above.
(See \cite[Lemma 3.1]{I1}.)

We expand our $Y_{lm}^n$ with respect to the basis $\{\ti T_e\}_e$.
Then the coefficient is given as follows using the relations
(\ref{Eq3+}), (\ref{Eq7}).

\begin{eqnarray}
\label{Eq8}
(Y_{lm}^n)^* \ti T_e &=& w^{-1/2}\phi_\Th(X_n^* X_l X_m) \ti T_e\nonumber\\
&=& w^{-1/2} V^* X_n^* X_l X_m V \ti T_e\nonumber\\
&=& \frac{w^{1/2}}{d_{\nu_1} d_{\nu_2}}
E_{M\otimes M^\o}( X_l X_m V \ti T_e V^* X_n^*)\nonumber\\
&=& \frac{1}{w^{1/2}d_{\nu_1} d_{\nu_2}}
X_l X_m \Th(\ti T_e)X_n^*
\end{eqnarray}

We represent $X_l$ graphically as in Fig.~\ref{Fg2},
where we follow the graphical convention of \cite[Section 3]{BEK1},
and $\{T_i\}_i$ is an orthonormal basis of 
$\Hom(\be_5, \be_1\be_3)$.  After this figure, we drop the symbols $T_i$,
$S_l^{\pm*}$, and the summation $\sum_{T_i}$ for simplicity.

\begin{figure}[htb]
\begin{center}
\unitlength 0.6mm
\begin{picture}(200,60)
\thinlines
\put(58,30){\makebox(0,0){$\displaystyle{X_l=\sum_{\be_3,\be_5}
\frac{\sqrt{d_{\la_1}d_{\la_2}}}{d_{\be_1}} W_{\be_5}(\sum_{T_i}}$}}
\put(143,30){\makebox(0,0){$\otimes\;j$}}
\put(193,30){\makebox(0,0){$)W_{\be_3}^*.$}}
\put(120,30){\arc{20}{0}{3.141593}}
\put(170,30){\arc{20}{0}{3.141593}}
\put(110,30){\circle*{2}}
\put(120,20){\circle*{2}}
\put(160,30){\circle*{2}}
\put(170,20){\circle*{2}}
\put(120,10){\line(0,1){10}}
\put(170,10){\line(0,1){10}}
\put(110,30){\line(1,1){20}}
\put(180,30){\line(-1,1){20}}
\put(110,50){\line(1,-1){8}}
\put(130,30){\line(-1,1){8}}
\put(160,30){\line(1,1){8}}
\put(180,50){\line(-1,-1){8}}
\put(120,6){\makebox(0,0){$\be_5$}}
\put(170,6){\makebox(0,0){$\be_5$}}
\put(110,54){\makebox(0,0){$\be_3$}}
\put(130,55){\makebox(0,0){$\a^+_{\la_1}$}}
\put(160,54){\makebox(0,0){$\be_3$}}
\put(180,55){\makebox(0,0){$\a^-_{\la_2}$}}
\put(110,20){\makebox(0,0){$\be_1$}}
\put(160,20){\makebox(0,0){$\be_1$}}
\put(106,32){\makebox(0,0){$S_l^{+*}$}}
\put(156,32){\makebox(0,0){$S_l^{-*}$}}
\put(126,15){\makebox(0,0){$T_i^*$}}
\put(176,15){\makebox(0,0){$T_i^*$}}
\end{picture}
\end{center}
\caption{A graphical expression for $X_l$}
\label{Fg2}
\end{figure}

We next have  a graphical expression for  $X_l X_m$ 
as in Fig.~\ref{Fg3},
where we have used a braiding-fusion equation for the
half-braiding.

\begin{figure}[htb]
\begin{center}
\unitlength 0.6mm
\begin{picture}(250,100)
\thinlines
\put(54,55){\makebox(0,0){$\displaystyle
{X_l X_m=\sum_{\be_3,\be'_3,\be_5}
\frac{\sqrt{d_{\la_1}d_{\la_2}d_{\mu_1}d_{\mu_2}}}
{d_{\be_1}d_{\be'_1}} W_{\be_5}(\sum_{T_i}}$}}
\put(172,55){\makebox(0,0){$\otimes\;j$}}
\put(238,55){\makebox(0,0){$)W_{\be'_3}^*.$}}
\put(135,30){\arc{30}{0}{3.141593}}
\put(150,50){\arc{20}{0}{3.141593}}
\put(205,30){\arc{30}{0}{3.141593}}
\put(220,50){\arc{20}{0}{3.141593}}
\put(135,15){\circle*{2}}
\put(120,30){\circle*{2}}
\put(150,40){\circle*{2}}
\put(120,70){\circle*{2}}
\put(205,15){\circle*{2}}
\put(190,30){\circle*{2}}
\put(220,40){\circle*{2}}
\put(190,70){\circle*{2}}
\put(135,5){\line(0,1){10}}
\put(205,5){\line(0,1){10}}
\put(120,30){\line(0,1){20}}
\put(190,30){\line(0,1){20}}
\put(150,30){\line(0,1){10}}
\put(220,30){\line(0,1){10}}
\put(120,70){\line(1,1){20}}
\put(120,50){\line(1,1){40}}
\put(120,90){\line(1,-1){8}}
\put(132,78){\line(1,-1){6}}
\put(142,68){\line(1,-1){18}}
\put(120,70){\line(1,-1){8}}
\put(132,58){\line(1,-1){8}}
\put(190,70){\line(1,1){8}}
\put(202,82){\line(1,1){8}}
\put(190,50){\line(1,1){8}}
\put(202,62){\line(1,1){6}}
\put(212,72){\line(1,1){18}}
\put(190,70){\line(1,-1){20}}
\put(190,90){\line(1,-1){40}}
\put(135,1){\makebox(0,0){$\be_5$}}
\put(205,1){\makebox(0,0){$\be_5$}}
\put(120,94){\makebox(0,0){$\be'_3$}}
\put(140,95){\makebox(0,0){$\a^+_{\mu_1}$}}
\put(160,95){\makebox(0,0){$\a^+_{\la_1}$}}
\put(190,94){\makebox(0,0){$\be'_3$}}
\put(210,95){\makebox(0,0){$\a^-_{\mu_2}$}}
\put(230,95){\makebox(0,0){$\a^-_{\la_2}$}}
\put(118,25){\makebox(0,0){$\be_1$}}
\put(188,25){\makebox(0,0){$\be_1$}}
\put(118,65){\makebox(0,0){$\be'_1$}}
\put(188,65){\makebox(0,0){$\be'_1$}}
\put(155,32){\makebox(0,0){$\be_3$}}
\put(225,32){\makebox(0,0){$\be_3$}}
\end{picture}
\end{center}
\caption{A graphical expression for $X_l X_m$}
\label{Fg3}
\end{figure}

Here we prepare two lemmas.

\begin{lemma}
\label{left}
For an intertwiner in $\Hom(\be_1\be_2,\be_3)\otimes
\Hom(\be_3,\be_1\be_2)$, the application of the left inverse
$\phi_{\be_1}$ is given as in Fig.~\ref{Fg11}.
\end{lemma}

\begin{proof}
Immediate by \cite[Lemma 3.1. (i)]{I1} and our graphical normalization
convention.
\end{proof}

\begin{figure}[htb]
\begin{center}
\unitlength 0.6mm
\begin{picture}(100,50)
\thinlines
\put(75,30){\arc{10}{3.141593}{6.283186}}
\put(75,20){\arc{10}{0}{3.141593}}
\put(30,20){\line(0,1){10}}
\put(70,20){\line(0,1){10}}
\put(80,20){\line(0,1){10}}
\put(20,10){\line(1,1){10}}
\put(30,30){\line(1,1){10}}
\put(80,30){\line(1,1){10}}
\put(40,10){\line(-1,1){10}}
\put(30,30){\line(-1,1){10}}
\put(90,10){\line(-1,1){10}}
\put(30,20){\circle*{2}}
\put(30,30){\circle*{2}}
\put(80,20){\circle*{2}}
\put(80,30){\circle*{2}}
\put(10,25){\makebox(0,0){$\phi_{\be_1}:$}}
\put(53,25){\makebox(0,0){$\mapsto\displaystyle\frac{1}{d_{\be_1}}$}}
\put(20,6){\makebox(0,0){$\be_1$}}
\put(40,6){\makebox(0,0){$\be_2$}}
\put(90,6){\makebox(0,0){$\be_2$}}
\put(20,44){\makebox(0,0){$\be_1$}}
\put(40,44){\makebox(0,0){$\be_2$}}
\put(90,44){\makebox(0,0){$\be_2$}}
\put(26,25){\makebox(0,0){$\be_3$}}
\put(66,30){\makebox(0,0){$\be_1$}}
\put(84,25){\makebox(0,0){$\be_3$}}
\end{picture}
\end{center}
\caption{A graphical expression for the left inverse}
\label{Fg11}
\end{figure}

\begin{figure}[htb]
\begin{center}
\unitlength 0.6mm
\begin{picture}(140,50)
\thinlines
\put(20,20){\line(0,1){10}}
\put(60,20){\line(0,1){10}}
\put(90,10){\line(0,1){30}}
\put(100,10){\line(0,1){30}}
\put(120,10){\line(0,1){30}}
\put(130,10){\line(0,1){30}}
\put(10,10){\line(1,1){10}}
\put(20,30){\line(1,1){10}}
\put(50,10){\line(1,1){10}}
\put(60,30){\line(1,1){10}}
\put(30,10){\line(-1,1){10}}
\put(20,30){\line(-1,1){10}}
\put(70,10){\line(-1,1){10}}
\put(60,30){\line(-1,1){10}}
\put(90,25){\line(1,0){10}}
\put(120,25){\line(1,0){10}}
\put(20,20){\circle*{2}}
\put(20,30){\circle*{2}}
\put(60,20){\circle*{2}}
\put(60,30){\circle*{2}}
\put(90,25){\circle*{2}}
\put(100,25){\circle*{2}}
\put(120,25){\circle*{2}}
\put(130,25){\circle*{2}}
\put(40,25){\makebox(0,0){$\otimes\;j$}}
\put(110,25){\makebox(0,0){$\otimes\;j$}}
\put(5,21){\makebox(0,0){$\displaystyle\sum_{\be_3}$}}
\put(80,21){\makebox(0,0){$=\displaystyle\sum_{\be_3}$}}
\put(10,6){\makebox(0,0){$\be_4$}}
\put(30,6){\makebox(0,0){$\be_5$}}
\put(50,6){\makebox(0,0){$\be_4$}}
\put(70,6){\makebox(0,0){$\be_5$}}
\put(90,6){\makebox(0,0){$\be_4$}}
\put(100,6){\makebox(0,0){$\be_5$}}
\put(120,6){\makebox(0,0){$\be_4$}}
\put(130,6){\makebox(0,0){$\be_5$}}
\put(10,44){\makebox(0,0){$\be_1$}}
\put(30,44){\makebox(0,0){$\be_2$}}
\put(50,44){\makebox(0,0){$\be_1$}}
\put(70,44){\makebox(0,0){$\be_2$}}
\put(90,44){\makebox(0,0){$\be_1$}}
\put(100,44){\makebox(0,0){$\be_2$}}
\put(120,44){\makebox(0,0){$\be_1$}}
\put(130,44){\makebox(0,0){$\be_2$}}
\put(16,25){\makebox(0,0){$\be_3$}}
\put(56,25){\makebox(0,0){$\be_3$}}
\put(95,30){\makebox(0,0){$\be'_3$}}
\put(125,30){\makebox(0,0){$\be'_3$}}
\end{picture}
\end{center}
\caption{A change of orthonormal bases}
\label{Fg12}
\end{figure}
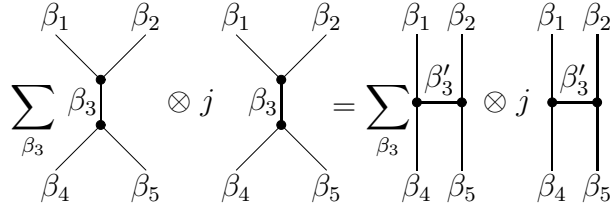

\begin{lemma}\label{summation}
For a change of bases, we have a graphical identity as in
Fig.~\ref{Fg12}, where we have summations over orthonormal bases
of (co)-isometries for small black circles.
\end{lemma}

\begin{proof}
The change of bases produces quantum $6j$-symbols, and their
unitarity gives the conclusion.
\end{proof}

\begin{figure}[htb]
\begin{center}
\unitlength 0.6mm
\begin{picture}(120,170)
\thinlines
\put(20,140){\arc{20}{3.141593}{6.283186}}
\put(20,50){\arc{20}{0}{3.141593}}
\put(35,40){\arc{30}{0}{3.141593}}
\put(35,10){\line(0,1){15}}
\put(10,50){\line(1,1){40}}
\put(30,50){\line(0,1){18}}
\put(30,90){\line(0,-1){18}}
\put(50,40){\line(0,1){30}}
\put(10,90){\line(0,1){30}}
\put(40,100){\line(0,1){10}}
\put(40,100){\line(1,-1){10}}
\put(40,100){\line(-1,-1){10}}
\put(10,140){\line(1,-1){30}}
\put(10,120){\line(1,1){8}}
\put(30,140){\line(-1,-1){8}}
\put(20,150){\line(0,1){10}}
\put(10,90){\line(2,-1){18}}
\put(32,79){\line(2,-1){2.66667}}
\put(50,70){\line(-2,1){11.3333}}
\put(10,140){\circle*{2}}
\put(20,150){\circle*{2}}
\put(40,100){\circle{2}}
\put(30,75){\circle*{2}}
\put(10,50){\circle*{2}}
\put(20,40){\circle*{2}}
\put(35,25){\circle*{2}}
\put(80,140){\arc{20}{3.141593}{6.283186}}
\put(80,50){\arc{20}{0}{3.141593}}
\put(95,40){\arc{30}{0}{3.141593}}
\put(95,10){\line(0,1){15}}
\put(70,50){\line(1,1){18}}
\put(92,71){\line(1,1){2.6667}}
\put(110,90){\line(-1,-1){11.3333}}
\put(90,50){\line(0,1){28}}
\put(90,90){\line(0,-1){8}}
\put(110,40){\line(0,1){30}}
\put(70,90){\line(0,1){30}}
\put(100,100){\line(0,1){10}}
\put(100,100){\line(1,-1){10}}
\put(100,100){\line(-1,-1){10}}
\put(70,140){\line(1,-1){8}}
\put(100,110){\line(-1,1){18}}
\put(70,120){\line(1,1){20}}
\put(80,150){\line(0,1){10}}
\put(70,90){\line(2,-1){40}}
\put(70,140){\circle*{2}}
\put(80,150){\circle*{2}}
\put(100,100){\circle{2}}
\put(90,75){\circle*{2}}
\put(70,50){\circle*{2}}
\put(80,40){\circle*{2}}
\put(95,25){\circle*{2}}
\put(60,90){\makebox(0,0){$\otimes\;j$}}
\put(35,6){\makebox(0,0){$\be_5$}}
\put(95,6){\makebox(0,0){$\be_5$}}
\put(19,30){\makebox(0,0){$\ti\be_3$}}
\put(79,30){\makebox(0,0){$\ti\be_3$}}
\put(52,30){\makebox(0,0){$\be'_3$}}
\put(112,30){\makebox(0,0){$\be'_3$}}
\put(8,44){\makebox(0,0){$\be_1$}}
\put(68,44){\makebox(0,0){$\be_1$}}
\put(35,50){\makebox(0,0){$\be'_1$}}
\put(95,50){\makebox(0,0){$\be'_1$}}
\put(19,67){\makebox(0,0){$\a^+_{\la_1}$}}
\put(79,67){\makebox(0,0){$\a^-_{\la_2}$}}
\put(29,97){\makebox(0,0){$\a^+_{\mu_1}$}}
\put(89,97){\makebox(0,0){$\a^-_{\mu_2}$}}
\put(46,107){\makebox(0,0){$\a^+_{\nu_1}$}}
\put(106,107){\makebox(0,0){$\a^-_{\nu_2}$}}
\put(35,142){\makebox(0,0){$\be'_3$}}
\put(95,142){\makebox(0,0){$\be'_3$}}
\put(7,145){\makebox(0,0){$\be''_1$}}
\put(67,145){\makebox(0,0){$\be''_1$}}
\put(20,164){\makebox(0,0){$\be_5$}}
\put(80,164){\makebox(0,0){$\be_5$}}
\end{picture}
\end{center}
\caption{A graphical expression for $X_l X_m\Th(\ti T_e)X_n^*$}
\label{Fg4}
\end{figure}
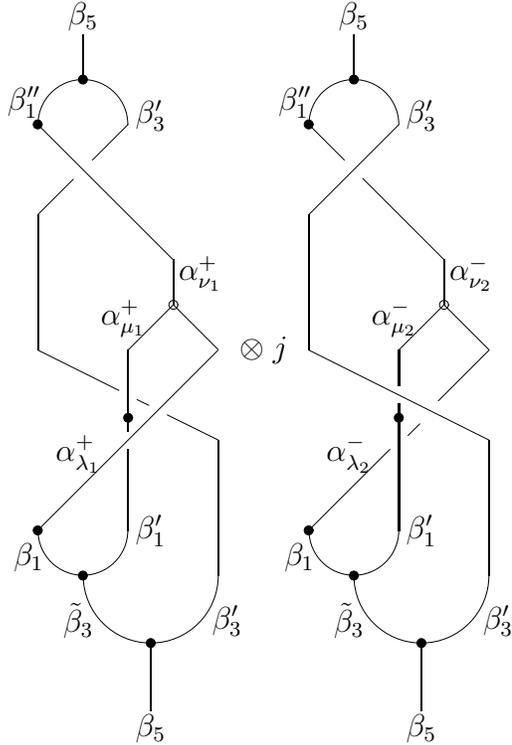

Then next we compute $X_l X_m\Th(\ti T_e)X_n^*$.
It is expressed as
\begin{eqnarray}
\label{Eq9}
&&X_l X_m\Th(\ti T_e)X_n^*\nonumber\\
&=& 
\frac{\sqrt{d_{\la_1}d_{\la_2}d_{\mu_1}d_{\mu_2}d_{\nu_1}d_{\nu_2}}}
{d_{\be_1}d_{\be'_1}d_{\be''_1}}\left(\frac{d_{\nu_1}d_{\nu_2}}
{d_{\la_1}d_{\la_2}d_{\mu_1}d_{\mu_2}}\right)^{1/4}\nonumber\\
&&\qquad\qquad\times
\sum_{\be'_3,\be_5,\ti\be_3} W_{\be_5}{\mathrm{
(graphical\ expression\ of\ Fig.~(\ref{Fg4}))}} W^*_{\be_5},
\end{eqnarray}
where small white circles represent intertwiners corresponding
to $T^1_{e_1}, T^2_{e^2}$ regarded as elements in $M$,
we have applied $\phi_\Th$ graphically using Lemma \ref{left},
changed the orthonormal bases in the space
$\Hom(\be_1\be'_1\be'_3,\be_5)$ using Lemma \ref{summation}
and thus we now have a summation
over $\ti\be_3$ rather than over $\be_3$.

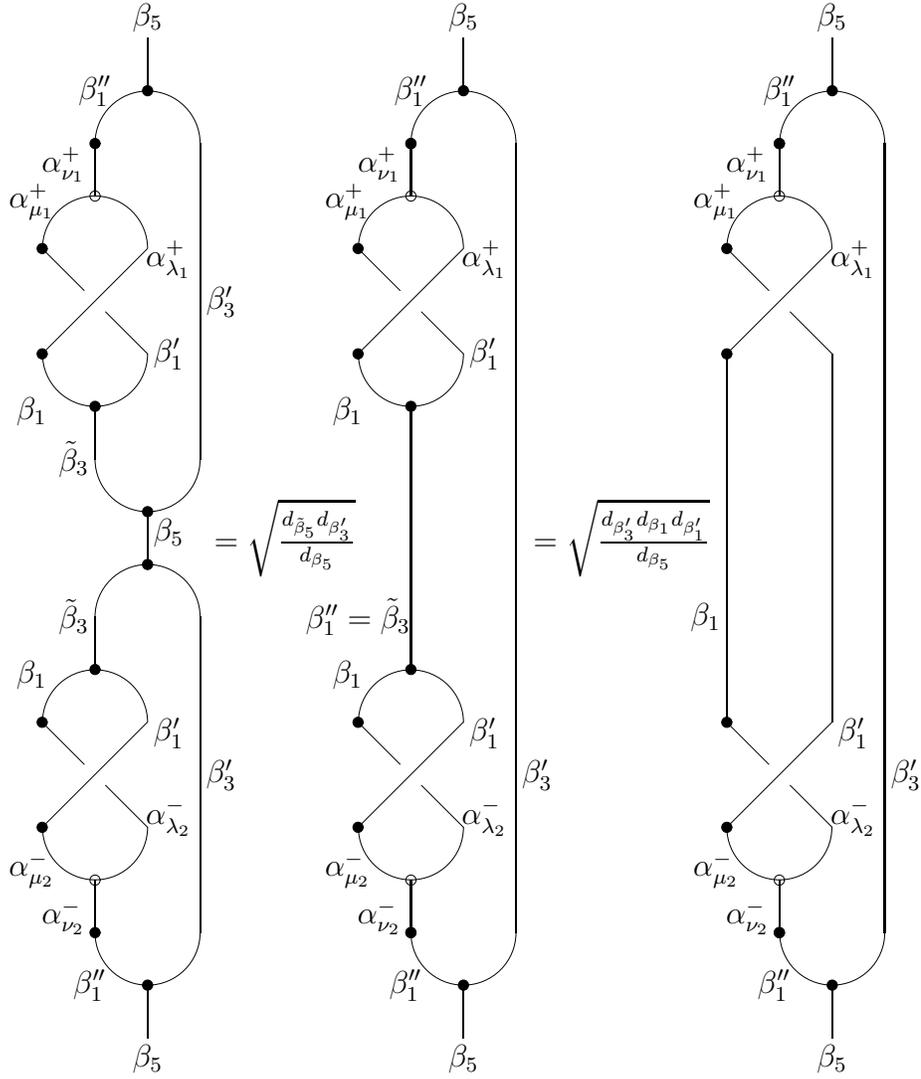
\begin{figure}[htb]
\begin{center}
\unitlength 0.7mm
\begin{picture}(180,210)
\thinlines
\put(20,160){\arc{20}{3.141593}{6.283186}}
\put(20,70){\arc{20}{3.141593}{6.283186}}
\put(30,180){\arc{20}{3.141593}{6.283186}}
\put(30,90){\arc{20}{3.141593}{6.283186}}
\put(90,180){\arc{20}{3.141593}{6.283186}}
\put(80,160){\arc{20}{3.141593}{6.283186}}
\put(80,70){\arc{20}{3.141593}{6.283186}}
\put(150,160){\arc{20}{3.141593}{6.283186}}
\put(160,180){\arc{20}{3.141593}{6.283186}}
\put(20,50){\arc{20}{0}{3.141593}}
\put(20,140){\arc{20}{0}{3.141593}}
\put(30,30){\arc{20}{0}{3.141593}}
\put(30,120){\arc{20}{0}{3.141593}}
\put(80,50){\arc{20}{0}{3.141593}}
\put(80,140){\arc{20}{0}{3.141593}}
\put(90,30){\arc{20}{0}{3.141593}}
\put(150,50){\arc{20}{0}{3.141593}}
\put(160,30){\arc{20}{0}{3.141593}}
\put(30,10){\line(0,1){10}}
\put(20,30){\line(0,1){10}}
\put(20,80){\line(0,1){10}}
\put(40,30){\line(0,1){60}}
\put(30,100){\line(0,1){10}}
\put(20,120){\line(0,1){10}}
\put(20,170){\line(0,1){10}}
\put(40,120){\line(0,1){60}}
\put(30,190){\line(0,1){10}}
\put(90,10){\line(0,1){10}}
\put(80,30){\line(0,1){10}}
\put(80,80){\line(0,1){50}}
\put(80,170){\line(0,1){10}}
\put(90,190){\line(0,1){10}}
\put(100,30){\line(0,1){150}}
\put(160,10){\line(0,1){10}}
\put(150,30){\line(0,1){10}}
\put(140,70){\line(0,1){70}}
\put(160,70){\line(0,1){70}}
\put(150,170){\line(0,1){10}}
\put(160,190){\line(0,1){10}}
\put(170,30){\line(0,1){150}}
\put(10,50){\line(1,1){20}}
\put(10,140){\line(1,1){20}}
\put(70,50){\line(1,1){20}}
\put(70,140){\line(1,1){20}}
\put(140,50){\line(1,1){20}}
\put(140,140){\line(1,1){20}}
\put(30,50){\line(-1,1){8}}
\put(10,70){\line(1,-1){8}}
\put(30,140){\line(-1,1){8}}
\put(10,160){\line(1,-1){8}}
\put(90,50){\line(-1,1){8}}
\put(70,70){\line(1,-1){8}}
\put(90,140){\line(-1,1){8}}
\put(70,160){\line(1,-1){8}}
\put(160,50){\line(-1,1){8}}
\put(140,70){\line(1,-1){8}}
\put(160,140){\line(-1,1){8}}
\put(140,160){\line(1,-1){8}}
\put(10,50){\circle*{2}}
\put(10,70){\circle*{2}}
\put(10,140){\circle*{2}}
\put(10,160){\circle*{2}}
\put(70,50){\circle*{2}}
\put(70,70){\circle*{2}}
\put(70,140){\circle*{2}}
\put(70,160){\circle*{2}}
\put(140,50){\circle*{2}}
\put(140,70){\circle*{2}}
\put(140,140){\circle*{2}}
\put(140,160){\circle*{2}}
\put(80,30){\circle*{2}}
\put(80,180){\circle*{2}}
\put(150,30){\circle*{2}}
\put(150,180){\circle*{2}}
\put(30,20){\circle*{2}}
\put(30,100){\circle*{2}}
\put(30,110){\circle*{2}}
\put(30,190){\circle*{2}}
\put(20,30){\circle*{2}}
\put(20,180){\circle*{2}}
\put(80,30){\circle{2}}
\put(80,180){\circle{2}}
\put(150,30){\circle{2}}
\put(150,180){\circle{2}}
\put(20,40){\circle{2}}
\put(20,80){\circle*{2}}
\put(20,130){\circle*{2}}
\put(20,170){\circle{2}}
\put(90,20){\circle*{2}}
\put(90,190){\circle*{2}}
\put(80,40){\circle{2}}
\put(80,80){\circle*{2}}
\put(80,130){\circle*{2}}
\put(80,170){\circle{2}}
\put(160,20){\circle*{2}}
\put(160,190){\circle*{2}}
\put(150,40){\circle{2}}
\put(150,170){\circle{2}}
\put(56,105){\makebox(0,0)
{$=\sqrt{\frac{d_{\ti\be_5} d_{\be'_3}}{d_{\be_5}}}$}}
\put(120,105){\makebox(0,0)
{$=\sqrt{\frac{d_{\be'_3} d_{\be_1} d_{\be'_1}}{d_{\be_5}}}$}}
\put(30,6){\makebox(0,0){$\be_5$}}
\put(90,6){\makebox(0,0){$\be_5$}}
\put(160,6){\makebox(0,0){$\be_5$}}
\put(19,20){\makebox(0,0){$\be''_1$}}
\put(79,20){\makebox(0,0){$\be''_1$}}
\put(149,20){\makebox(0,0){$\be''_1$}}
\put(14,33){\makebox(0,0){$\a^-_{\nu_2}$}}
\put(74,33){\makebox(0,0){$\a^-_{\nu_2}$}}
\put(144,33){\makebox(0,0){$\a^-_{\nu_2}$}}
\put(8,42){\makebox(0,0){$\a^-_{\mu_2}$}}
\put(68,42){\makebox(0,0){$\a^-_{\mu_2}$}}
\put(138,42){\makebox(0,0){$\a^-_{\mu_2}$}}
\put(34,52){\makebox(0,0){$\a^-_{\la_2}$}}
\put(94,52){\makebox(0,0){$\a^-_{\la_2}$}}
\put(164,52){\makebox(0,0){$\a^-_{\la_2}$}}
\put(34,68){\makebox(0,0){$\be'_1$}}
\put(94,68){\makebox(0,0){$\be'_1$}}
\put(164,68){\makebox(0,0){$\be'_1$}}
\put(8,79){\makebox(0,0){$\be_1$}}
\put(68,79){\makebox(0,0){$\be_1$}}
\put(136,90){\makebox(0,0){$\be_1$}}
\put(16,90){\makebox(0,0){$\ti\be_3$}}
\put(16,120){\makebox(0,0){$\ti\be_3$}}
\put(70,90){\makebox(0,0){$\be''_1=\ti\be_3$}}
\put(34,106){\makebox(0,0){$\be_5$}}
\put(44,60){\makebox(0,0){$\be'_3$}}
\put(44,150){\makebox(0,0){$\be'_3$}}
\put(104,60){\makebox(0,0){$\be'_3$}}
\put(174,60){\makebox(0,0){$\be'_3$}}
\put(8,129){\makebox(0,0){$\be_1$}}
\put(68,129){\makebox(0,0){$\be_1$}}
\put(34,140){\makebox(0,0){$\be'_1$}}
\put(94,140){\makebox(0,0){$\be'_1$}}
\put(34,158){\makebox(0,0){$\a^+_{\la_1}$}}
\put(94,158){\makebox(0,0){$\a^+_{\la_1}$}}
\put(164,158){\makebox(0,0){$\a^+_{\la_1}$}}
\put(8,169){\makebox(0,0){$\a^+_{\mu_1}$}}
\put(68,169){\makebox(0,0){$\a^+_{\mu_1}$}}
\put(138,169){\makebox(0,0){$\a^+_{\mu_1}$}}
\put(14,176){\makebox(0,0){$\a^+_{\nu_1}$}}
\put(74,176){\makebox(0,0){$\a^+_{\nu_1}$}}
\put(144,176){\makebox(0,0){$\a^+_{\nu_1}$}}
\put(20,190){\makebox(0,0){$\be''_1$}}
\put(80,190){\makebox(0,0){$\be''_1$}}
\put(150,190){\makebox(0,0){$\be''_1$}}
\put(30,204){\makebox(0,0){$\be_5$}}
\put(90,204){\makebox(0,0){$\be_5$}}
\put(160,204){\makebox(0,0){$\be_5$}}
\end{picture}
\end{center}
\caption{The value of Fig.~\ref{Fg4}}
\label{Fg5}
\end{figure}

Then the complex number value represented by
Fig.~\ref{Fg4} can be
computed as in Fig.~\ref{Fg5}, where we have used the 
braiding-fusion equation for a half-braiding twice.

Here we have the following lemma.

\begin{figure}[htb]
\begin{center}
\unitlength 0.5mm
\begin{picture}(120,60)
\thinlines
\put(10,50){\line(1,-1){40}}
\put(70,50){\line(1,-1){40}}
\put(10,10){\line(1,1){18}}
\put(50,50){\line(-1,-1){18}}
\put(70,10){\line(1,1){18}}
\put(110,50){\line(-1,-1){18}}
\put(20,20){\circle*{2}}
\put(20,40){\circle*{2}}
\put(80,40){\circle*{2}}
\put(100,40){\circle*{2}}
\put(10,4){\makebox(0,0){$\be'$}}
\put(50,4){\makebox(0,0){$\be$}}
\put(70,4){\makebox(0,0){$\be'$}}
\put(110,4){\makebox(0,0){$\be$}}
\put(10,56){\makebox(0,0){$\a^+_\la$}}
\put(50,56){\makebox(0,0){$\a^+_\mu$}}
\put(70,56){\makebox(0,0){$\a^+_\la$}}
\put(110,56){\makebox(0,0){$\a^+_\mu$}}
\put(12,40){\makebox(0,0){$T^*$}}
\put(12,20){\makebox(0,0){$S^*$}}
\put(72,40){\makebox(0,0){$T^*$}}
\put(108,40){\makebox(0,0){$S^*$}}
\put(60,30){\makebox(0,0){$=$}}
\end{picture}
\end{center}
\caption{A naturality equation}
\label{Fg8}
\end{figure}
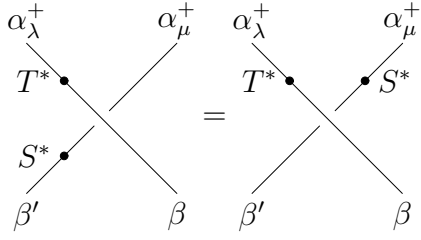

\begin{lemma}
\label{nat}
Let $\be,\be'$ be ambichiral and choose isometries
$T\in\Hom(\be,\a^+_\la)$, $S\in\Hom(\be',\a^+_\mu)$.
Then we have the identity as in Fig.~\ref{Fg8}.
\end{lemma}

\begin{proof}
We compute the both hand sides 
by the definitions of the half and
the relative braidings in \cite[(10)]{BEK3} and
\cite[Subsection 3.3]{BE}, respectively, and
then we get $\be'(T^*)S^* \e^+(\la,\mu) TT^*$, where
we have used $\e^+(\la,\mu)\a_\la^+(SS^*)=
SS^*\e^+(\la,\mu)$, which follows from the arguments and
the figure in \cite[page 377]{X1}.
(The chiral locality is not used
in the argument in \cite[page 377]{X1}.)
\end{proof}

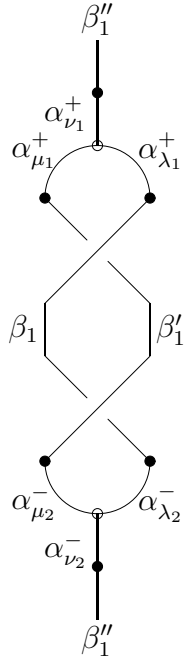
\begin{figure}[htb]
\begin{center}
\unitlength 0.7mm
\begin{picture}(40,130)
\thinlines
\put(20,90){\arc{20}{3.141593}{6.283186}}
\put(20,40){\arc{20}{0}{3.141593}}
\put(20,10){\line(0,1){20}}
\put(20,100){\line(0,1){20}}
\put(10,60){\line(0,1){10}}
\put(30,60){\line(0,1){10}}
\put(10,40){\line(1,1){20}}
\put(10,70){\line(1,1){20}}
\put(30,40){\line(-1,1){8}}
\put(10,60){\line(1,-1){8}}
\put(30,70){\line(-1,1){8}}
\put(10,90){\line(1,-1){8}}
\put(10,40){\circle*{2}}
\put(30,40){\circle*{2}}
\put(30,90){\circle*{2}}
\put(10,90){\circle*{2}}
\put(20,20){\circle*{2}}
\put(20,30){\circle{2}}
\put(20,100){\circle{2}}
\put(20,110){\circle*{2}}
\put(20,6){\makebox(0,0){$\be''_1$}}
\put(20,124){\makebox(0,0){$\be''_1$}}
\put(14,23){\makebox(0,0){$\a^-_{\nu_2}$}}
\put(8,32){\makebox(0,0){$\a^-_{\mu_2}$}}
\put(32,32){\makebox(0,0){$\a^-_{\la_2}$}}
\put(34,65){\makebox(0,0){$\be'_1$}}
\put(6,65){\makebox(0,0){$\be_1$}}
\put(32,99){\makebox(0,0){$\a^+_{\la_1}$}}
\put(8,99){\makebox(0,0){$\a^+_{\mu_1}$}}
\put(14,106){\makebox(0,0){$\a^+_{\nu_1}$}}
\end{picture}
\end{center}
\caption{The remaining intertwiner}
\label{Fg6}
\end{figure}

Then the value $(Y_{lm}^n)^* \ti T_e$ is computed with 
the coefficients in the equations (\ref{Eq8}), (\ref{Eq9}),
and Fig.~\ref{Fg5}.  The coefficient is now
\begin{eqnarray}
\label{Eq10}
&&\frac{w^{1/2}}{d_{\nu_1} d_{\nu_2}}
\frac{\sqrt{d_{\la_1}d_{\la_2}d_{\mu_1}d_{\mu_2}d_{\nu_1}d_{\nu_2}}}
{d_{\be_1}d_{\be'_1}d_{\be''_1}}\left(\frac{d_{\nu_1}d_{\nu_2}}
{d_{\la_1}d_{\la_2}d_{\mu_1}d_{\mu_2}}\right)^{1/4}\nonumber\\
&&\qquad\qquad\times
\sqrt{\frac{d_{\be'_3} d_{\be_1} d_{\be'_1}}{d_{\be_5}}}
\sqrt{\frac{d_{\be''_1} d_{\be_5}}{d_{\be'_3}}} d_{\be''_1}
\nonumber\\
&=&w^{-1/2}
\sqrt{\frac{d_{\la_1}d_{\la_2}d_{\mu_1}d_{\mu_2}}
{d_{\nu_1}d_{\nu_2}}}
\sqrt{\frac{d_{\be''_1}}{d_{\be_1}d_{\be'_1}}}
\left(\frac{d_{\nu_1}d_{\nu_2}}
{d_{\la_1}d_{\la_2}d_{\mu_1}d_{\mu_2}}\right)^{1/4}
\end{eqnarray}
and this is multiplied with the intertwiner in Fig.~\ref{Fg6},
where the two crossings of the two wires labeled with
$\be_1,\be'_1$ represent the ``ambichiral braiding''
studied in \cite[Subsection 3.3]{BE}.

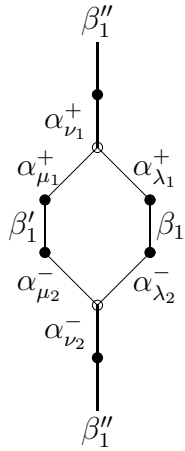
\begin{figure}[htb]
\begin{center}
\unitlength 0.7mm
\begin{picture}(40,90)
\thinlines
\put(20,10){\line(0,1){20}}
\put(20,60){\line(0,1){20}}
\put(10,40){\line(0,1){10}}
\put(30,40){\line(0,1){10}}
\put(20,30){\line(1,1){10}}
\put(20,60){\line(1,-1){10}}
\put(20,30){\line(-1,1){10}}
\put(20,60){\line(-1,-1){10}}
\put(20,20){\circle*{2}}
\put(20,30){\circle{2}}
\put(20,60){\circle{2}}
\put(20,70){\circle*{2}}
\put(10,40){\circle*{2}}
\put(10,50){\circle*{2}}
\put(30,40){\circle*{2}}
\put(30,50){\circle*{2}}
\put(20,6){\makebox(0,0){$\be''_1$}}
\put(20,84){\makebox(0,0){$\be''_1$}}
\put(14,25){\makebox(0,0){$\a^-_{\nu_2}$}}
\put(9,34){\makebox(0,0){$\a^-_{\mu_2}$}}
\put(31,34){\makebox(0,0){$\a^-_{\la_2}$}}
\put(34,45){\makebox(0,0){$\be_1$}}
\put(6,45){\makebox(0,0){$\be'_1$}}
\put(31,56){\makebox(0,0){$\a^+_{\la_1}$}}
\put(9,56){\makebox(0,0){$\a^+_{\mu_1}$}}
\put(14,65){\makebox(0,0){$\a^+_{\nu_1}$}}
\end{picture}
\end{center}
\caption{The new form of the remaining intertwiner}
\label{Fg7}
\end{figure}

Then the monodromy of $\be'_1$ and $\be_1$ in Fig.~\ref{Fg6} 
acts on $\Hom(\be'_1\be_1,\be''_1)$ as a scalar
arising from ``conformal dimensions'' of
$\be_1,\be'_1,\be''_1$ in the ambichiral system.  (See
\cite[Figure 8.30]{EK}.)  So
up to this scalar, we have Fig.~\ref{Fg7}.
Since the fourth root in (\ref{Eq10}) comes from our
normalization for the graphical expression (see
\cite[Figures 7, 9]{BEK1}) and
we can absorb the above scalar arising from the
conformal dimensions by changing the bases $\{\ti T_e\}_e$,
our coefficient multiplied with the number represented by
Fig.~\ref{Fg6} now coincides with Rehren's coefficient
computed as in (\ref{Rc}).
(Actually, $\la_j$ and $\mu_j$ are  interchanged and also
$\a^+$ and $\a^-$ are interchanged, but these are just matters
of convention.)

Now with \cite[Corollary 3.10]{BEK3}, we have proved the following
theorem.

\begin{theorem}
The generalized Longo-Rehren subfactor arising from 
$\a^\pm$-induction with a non-degenerate braiding on $\NXN$ is
isomorphic to the dual of the Longo-Rehren subfactor
arising from $\MXM$.
\end{theorem}

At the end of \cite{R2}, Rehren asks for an Izumi type
description \cite{I1} of irreducible endomorphisms of $P$ arising
from the generalized Longo-Rehren subfactor
$N\otimes N^\opp\subset P$ and
in particular, he asks whether a braiding exists or not
on this system of endomorphisms of $P$.   The above
theorem in particular shows that the system of endomorphisms of
$P$ is isomorphic to the direct product system of
$\MXM$ and $\MXM^\opp$ and thus
we solve these problems and the answer
to the second question is negative, since this
system does not have a braiding in general and it
can be even non-commutative.  (Note that
\cite[Corollary 6.9]{BEK1} gives a criterion for such
non-commutativity.)

\begin{remark}\label{tw}
{\rm If $N=M$ in the above setting, our result implies
\cite[Proposition 7.3]{I1}, of course, but
a remark on \cite[page 171]{I1} gives a 
``twisted Longo-Rehren subfactor'' rather than the
usual Longo-Rehren subfactor.   This is due to the
monodoromy operator similar to the one in Fig.~\ref{Fg6},
but as pointed by
Rehren, one can always eliminate such a twist and then
the  ``twisted Longo-Rehren subfactor'' is actually
isomorphic to the Longo-Rehren subfactor.  (See
``Added in proof'' of \cite{I1} on this point.)
We also had a similar twist in our results here,
originally, but we have eliminated it thanks to this
remark of Rehren.
}\end{remark}

In the above setting, we can also set
$N_1=N$, $N_2=M$, $\Delta_1=\NXN$, $\Delta_2={}_M{\mathcal{X}^0}_M$,
$\a^1_\la=\a^+_\la$, $\a^2_\tau=\tau$ in the construction of the generalized
Longo-Rehren subfactor.  Then the resulting subfactor
$M\otimes N^\o\subset R$ has a dual canonical endomorphism
$\bigoplus_{\la\in\NXN,\tau\in {}_M{\mathcal{X}^0}_M}
b^+_{\tau,\la} \la\otimes \tau^\o$,
where $b^+_{\tau,\la}=\dim\Hom(\a^+_\la,\tau)$
is the {\sl chiral branching coefficient} as in \cite[Subsection 3.2]{BEK2}.
Now using the results in \cite[Section 4]{BEK3} and arguments almost identical
to the above, we can prove the following theorem.

\begin{theorem}
The generalized Longo-Rehren subfactor $M\otimes N^\o\subset R$ arising from 
$\a^+$-induction as above with a non-degenerate braiding on $\NXN$ is
isomorphic to the dual of the Longo-Rehren subfactor
arising from ${}_M{\mathcal{X}^+}_M$.
\end{theorem}

\section{Nets of subfactors on $S^1$}

In this section, we study multi-interval subfactors
for completely rational nets of subfactors,  which
generalizes the study in \cite{KLM}.

Let $\{M(I)\}_{I\subset S^1}$ be a completely
rational net of factors of $S^1$ in the sense of
\cite{KLM}, where an ``interval''
$I$ is a non-empty, non-dense connected open subset
of $S^1$.  (That is, we assume isotony, conformal invariance,
positivity of the energy, locality, existence of the vacuum,
irreducibility, the split property,
strong additivity, and finiteness of the $\mu$-index.
See \cite{GL,KLM} for the detailed definitions.)
We also suppose to have a conformal subnet
$\{N(I)\}_{I\subset S^1}$  of $\{M(I)\}_{I\subset S^1}$
with finite index as in \cite{L4}.  The main result
in \cite{L4} says that the subnet $\{N(I)\}_{I\subset S^1}$
is also completely rational.

Let $E=I_1\cup I_3$ be a union of two intervals $I_1$, $I_3$
such that $\bar I_1 \cap \bar I_3 =\emptyset$.  Label the
interiors of the two connected components of $S\setminus E$
as $I_2$, $I_4$ so that $I_1,I_2, I_3,I_4$  appear
on the circle in a counterclockwise order.
We set $N_j=N(I_j)$, $M_j=M(I_j)$, for $j=1,2,3,4$.  (This
numbering should not be confused with the basic
construction.)  We also set $N=N_1, M=M_1$.

We have a finite system of mutually inequivalent irreducible 
DHR endomorphisms $\{\la\}$ for the net $\{N(I)\}$ by
complete rationality.
We may and do regard this as a braided system of endomorphisms
of $N=N_1$.  By \cite[Corollary 37]{KLM}, this braiding
is non-degenerate.  We write $\NXN$ for this system.
As in \cite{BEK1}, we can apply $\a^\pm$-induction to
get systems $\MXM,  {}_M{\mathcal{X}^+}_M,
{}_M{\mathcal{X}^-}_M,  {}_M{\mathcal{X}^0}_M$ of
irreducible endomorphisms of $M$.  That is, they are the systems of
irreducible endomorphisms of $M$ arising from $\a^\pm$-induction,
$\a^+$-induction, $\a^-$-induction, and the ``ambichiral''
system, respectively.
Since the braiding on $\NXN$ is non-degenerate,
\cite[Theorem 5.10]{BEK1} and
\cite[Proposition 5.1]{BE} imply that the ambichiral system
${}_M{\mathcal{X}^0}_M$ is given by the irreducible
DHR endomorphisms of the net $\{M(I)\}$.
By the inclusions
${}_M{\mathcal{X}^0}_M\subset {}_M{\mathcal{X}^\pm}_M\subset
\MXM$ and the Galois correspondence of \cite[Theorem 2.5]{I1}
(or by the characterization of the Longo-Rehren subfactor in
\cite[Appendix A]{KLM}), we have inclusions of the
corresponding Longo-Rehren subfactors $M\otimes M^\o\subset R$,
$M\otimes M^\o\subset R^\pm$, $M\otimes M^\o\subset R^0$ with
$R^0\subset R^\pm \subset R$.  We study these Longo-Rehren
subfactors in connection to the results in Section 2.

As in \cite{KLM}, we make identification of $S^1$ with
$\R\cup\{\infty\}$, and
as in \cite[Proposition 36]{KLM}, we may and do assume that
$I_1=(-b,-a),I_3=(a,b)$, with $0<a<b$.
Take a DHR endomorphism $\la$ localized in $I_1$ for the 
net $\{N(I)\}$. Let $P=M(\ti I)$, where $\ti I=(-\infty,0)$.
Let $J$ be the modular conjugation for $P$
with respect to the vacuum vector.
We consider endomorphisms of the $C^*$-algebras
$\overline{\bigcup_{I\subset (-\infty,\infty)} M(I)}$ and
$\overline{\bigcup_{I\subset (-\infty,\infty)} N(I)}$.  The
canonical endomorphism $\ga$ and the dual canonical endomorphism
$\th$ are regarded as endomorphisms of these $C^*$-algebras.
We regard $\a^+_\la$ as an endomorphism of the former $C^*$-algebra
as in \cite{LR},
and then it is not localized in $I_1$ any more, but
it is localized in $(-\infty,-a)$
by \cite[Proposition 3.9]{LR}.  We study an irreducible decomposition
of $\a^+_\la$ as an endomorphism of $M_1$ and choose 
$\be$ appearing in such an irreducible decomposition of $\a^+_\la$
regarded as an endomorphism of $M_1$.  That is, we choose an
isometry $W\in M_1$ with $W^*W\in\a_\la^+(M)'\cap  M$,
$\be(x)=W^* \a^+_\la(x) W$.  Using
this same formula, we can regard $\be$ as endomorphism of the
$C^*$-algebra $\overline{\bigcup_{I\subset (-\infty,\infty)} M(I)}$.
We next regard $\be$ as an endomorphism of $P$ and let $V_\be$ be
the isometry standard implementation of $\be\in\End(P)$ as
in \cite[Appendix]{GL}.  We now set
$\bar\be=J\be J$.  Then for any $X\in P\vee P'$, we have
$\be\bar\be(X)V_\be=V_\be X$ as in the proof of 
\cite[Proposition 36]{KLM} since $JV_\be J=V_\be$.  By strong
additivity, we have this for all local operators $X$.
Since $\la,\bar\la=J\la J, \be,\bar\be$ are 
localized in $(-\infty,a)$, $(a,\infty)$, $I_1$, $I_3$,
respectively, we know that $V_\be\in (M_2\vee N_4)'$.
Consider the subfactor $M_1\vee M_3\subset (M_2\vee N_4)'$.
By Frobenius reciprocity \cite{I0}, we know that the
dual canonical endomorphism for the subfactor
$M_1\vee M_3\subset (M_2\vee N_4)'$ contains
$\be\otimes\be^\o$, where $M_3=J M_1 J$ is now regarded
as $M_1^\o$ and $M_1\vee M_3$ is regarded as
$M_1\otimes M_1^\o$, for all $\be\in {}_M{\mathcal{X}^+}_M$.
We now compute the index of the subfactor 
$M_1\vee M_3\subset (M_2\vee N_4)'$ in two ways.
On one hand, it has an intermediate subfactor $(M_2\vee M_4)'$
and the index for $M_1\vee M_3\subset (M_2\vee M_4)'$ is
the global index of the ambichiral system by
\cite[Theorem 33]{KLM}.
The index of
$(M_2\vee M_4)'\subset (M_2\vee N_4)'$ is simply that of
the net $\{N(I)\subset M(I)\}$ of subfactors.
We also have
$$\frac{w_+}{w_0}=\frac{w}{w_+}=\sum_{\la\in\NXN} d_\la Z_{\la 0}
=d_\th=[M(I):N(I)],$$
where $w, w_+, w_0$ are the global indices of
$\MXM, {}_M{\mathcal{X}^+}_M, {}_M{\mathcal{X}^0}_M$,
respectively, by \cite[Theorem 4.2, Proposition 3.1]{BEK2},
\cite[Theorem 3.3 (1)]{X1}.  (Here we have used the chiral locality
condition arising from the locality of the net $\{M(I)\}$.
Without the chiral locality, the results in this section would
not hold in general.)
These imply that
\begin{equation}\label{g-ind}
[(M_2\vee N_4)':M_1\vee M_3]=w_+.
\end{equation}

On the other hand, the dual
canonical endomorphism for the subfactor
$M_1\vee M_3\subset (M_2\vee N_4)'$ contains
$\bigoplus_{\be \in {}_M{\mathcal{X}^+}_M} \be\otimes\be^\o$ from
the above considerations since each $\be$ is irreducible as an
endomorphism of $M$,
thus the index value is at least $\sum_{\be \in {}_M{\mathcal{X}^+}_M}
d_\be^2=w_+$.  Together with (\ref{g-ind}), we know that
the dual canonical endomorphism is indeed equal to
$\bigoplus_{\be \in {}_M{\mathcal{X}^+}_M} \be\otimes\be^\o$.

Put $R_\be=\sqrt{d_\be}V_\be\in (M_2\vee N_4)'$.  As in the proof of
\cite[Proposition 36]{KLM}, we now conclude that
the subfactor $M_1\vee M_3\subset (M_2\vee N_4)'$
is isomorphic to the Longo-Rehren subfactor
$M\otimes M^\o\subset R^+$.
Similarly, we know that the subfactor
$M_1\vee M_3\subset (N_2\vee M_4)'$ is isomorphic to
the Longo-Rehren subfactor $M\otimes M^\o\subset R^-$.
These two isomorphisms are compatible on $(M_2\vee M_4)'$
and they give an isomorphism of $M_1\vee M_3\subset
(M_2\vee M_4)'$ to the Longo-Rehren subfactor
$M\otimes M^\o\subset R^0$.
We finally look at the inclusions
\begin{eqnarray}
\label{comm}
\begin{array}{ccccc}
M\otimes M^\o &\subset & R^0 & \subset & R^+ \nonumber\\
&& \cap && \cap \nonumber\\
&& R^- & \subset & R.
\end{array}
\end{eqnarray}
The right square is a commuting square by \cite[Lemma 1]{L4}
and thus $R$ is generated by $R^+$ and $R^-$. 
(Or \cite[Theorem 5.10]{BEK1}
and \cite[Proposition 2.4, Theorem 2.5]{I1} also give this
generating property.)
It means that the above isomorphisms give the following
theorem.

\begin{theorem}
Under the above setting, 
the following system of algebras arising from four intervals on
the circle is isomorphic to 
the system of algebras (\ref{comm}) arising as Longo-Rehren subfactors.
$$\begin{array}{ccccc}
M_1\vee M_3&\subset & (M_2\vee M_4)' & \subset & (M_2\vee N_4)' \\
&& \cap && \cap \\
&& (N_2\vee M_4)' & \subset & (N_2\vee N_4)'.
\end{array}$$
\end{theorem}

\begin{remark}\label{rem2}{\rm
Passing to the commutant, we also conclude that
the subfactor $N_1\vee N_3\subset(M_2\vee M_4)'$
is isomorphic to the dual of $M\otimes M^\o\subset R$
and thus isomorphic to the generalized Longo-Rehren subfactor
arising from the $\a^\pm$-induction studied in Section 2.
In the example of the
conformal inclusion $SU(2)_{10}\subset Spin(5)_1$ in
\cite[Section 4.1]{X1},
this fact was first noticed by Rehren and it
can be proved also in general directly by computing the corresponding
$Q$-system.
}\end{remark}

\begin{footnotesize}
\vspace{0.5cm}
\noindent{\it Acknowledgment.}
The author thanks K.-H. Rehren for his remarks 
mentioned in Remarks \ref{tw}, \ref{rem2} and detailed
comments on a preliminary version of this paper.  We also
thank F. Xu for his comments on the preliminary version.
We gratefully acknowledge the financial supports of
Grant-in-Aid for Scientific Research, Ministry of
Education and Science (Japan), Japan-Britain joint
research project (2000 April--2002 March) of
Japan Society for the Promotion of Science,
Mathematical Sciences Research Institute (Berkeley),
the Mitsubishi Foundation and
University of Tokyo.  A part of this work was
carried out at Mathematical Sciences Research Institute, Berkeley,
and Universit\`a di Roma ``Tor Vergata'' and we thank them
for their hospitality.

\end{footnotesize}
\end{document}